\documentclass{amsart}
\usepackage{amsfonts}
\usepackage{graphicx}
\usepackage{amscd}

\newtheorem{theorem}{Theorem}
\theoremstyle{plain}
\newtheorem{acknowledgement}{Acknowledgement}

\newtheorem{lemma}{Lemma}

\newtheorem{remark}{Remark}

\numberwithin{equation}{section}

\input{tcilatex}

\begin{document}
\title[\v{C}eby\v{s}ev's Inequality for Weighted Means]{On the \v{C}eby\v{s}%
ev's Inequality for Weighted Means }
\author{S.S. Dragomir}
\address{School of Computer Science and Mathematics\\
Victoria University of Technology\\
PO Box 14428\\
Melbourne City MC 8001\\
Victoria, Australia.}
\email{sever.dragomir@vu.edu.au}
\urladdr{http://rgmia.vu.edu.au/SSDragomirWeb.html}
\date{March 03, 2003.}
\subjclass{Primary 26D15; Secondary 26D10.}
\keywords{\v{C}eby\v{s}ev's Inequality.}

\begin{abstract}
Some new sufficient conditions for the weighted \v{C}eby\v{s}ev's inequality
for real numbers to hold are provided.
\end{abstract}

\maketitle

\section{Introduction}

Consider the real sequences $\left( n-\text{tuples}\right) $ $\mathbf{a}%
=\left( a_{1},\dots ,a_{n}\right) ,$ $\mathbf{b}=\left( b_{1},\dots
,b_{n}\right) $ and the nonnegative sequence $\mathbf{p}=\left( p_{1},\dots
,p_{n}\right) $ with $P_{n}:=\sum_{i=1}^{n}p_{i}>0.$ Define the \textit{%
weighted \v{C}eby\v{s}ev's functional} 
\begin{equation}
T_{n}\left( \mathbf{p};\mathbf{a},\mathbf{b}\right) :=\frac{1}{P_{n}}%
\sum_{i=1}^{n}p_{i}a_{i}b_{i}-\frac{1}{P_{n}}\sum_{i=1}^{n}p_{i}a_{i}\cdot 
\frac{1}{P_{n}}\sum_{i=1}^{n}p_{i}b_{i}.  \label{a.1}
\end{equation}

In 1882 -- 1883, \v{C}eby\v{s}ev \cite{CE1} and \cite{CE2} proved that if $%
\mathbf{a}$ and $\mathbf{b}$ are monotonic in the same (opposite) sense,
then 
\begin{equation}
T_{n}\left( \mathbf{p};\mathbf{a},\mathbf{b}\right) \geq \left( \leq \right)
0.  \label{a.2}
\end{equation}

In the special case $\mathbf{p}=\mathbf{a}\geq \mathbf{0}$, it appears that
the inequality (\ref{a.2}) has been obtained by Laplace long before \v{C}eby%
\v{s}ev (see for example \cite[p. 240]{MPF1}).

The inequality (\ref{a.2}) was mentioned by Hardy, Littlewood and P\'{o}lya
in their book \cite{HLP} in 1934 in the more general setting of synchronous
sequences, i.e., if $\mathbf{a},$ $\mathbf{b}$ are synchronous
(asynchronous), this means that 
\begin{equation}
\left( a_{i}-a_{j}\right) \left( b_{i}-b_{j}\right) \geq \left( \leq \right)
0\text{ for any }i,j\in \left\{ 1,\dots ,n\right\} ,  \label{a.3}
\end{equation}%
then (\ref{a.2}) holds true as well.

A relaxation of the synchronicity condition was provided by M. Biernacki in
1951, \cite{BIE}, which showed that, if $\mathbf{a},$ $\mathbf{b}$ are
monotonic in mean in the same sense, i.e., for $P_{k}:=\sum_{i=1}^{k}p_{i},$ 
$k=1,\dots ,n-1;$%
\begin{equation}
\frac{1}{P_{k}}\sum_{i=1}^{k}p_{i}a_{i}\leq \left( \geq \right) \frac{1}{%
P_{k+1}}\sum_{i=1}^{k+1}p_{i}a_{i},\;\;k\in \left\{ 1,\dots ,n-1\right\} 
\label{a.4}
\end{equation}%
and 
\begin{equation}
\frac{1}{P_{k}}\sum_{i=1}^{k}p_{i}b_{i}\leq \left( \geq \right) \frac{1}{%
P_{k+1}}\sum_{i=1}^{k+1}p_{i}b_{i},\;\;k\in \left\{ 1,\dots ,n-1\right\} ,
\label{a.5}
\end{equation}%
then (\ref{a.2}) holds with \textquotedblleft $\ \geq $ \ \textquotedblright
. If if $\mathbf{a},$ $\mathbf{b}$ are monotonic in mean in the opposite
sense then (\ref{a.2}) holds with \textquotedblleft $\ \leq $ \
\textquotedblright .

In 1989, Dragomir and Pe\v{c}ari\'{c} \cite{5b} proved the following
refinement of \v{C}eby\v{s}ev's inequality for synchronous sequences. If $%
\mathbf{a},$ $\mathbf{b}$ are synchronous and by $\left\vert \mathbf{a}%
\right\vert $ we denote the $n-$tuple $\left( \left\vert a_{1}\right\vert
,\dots ,\left\vert a_{n}\right\vert \right) ,$ then 
\begin{equation}
T_{n}\left( \mathbf{p};\mathbf{a},\mathbf{b}\right) \geq \max \left\{
\left\vert T_{n}\left( \mathbf{p};\left\vert \mathbf{a}\right\vert ,\mathbf{b%
}\right) \right\vert ,\left\vert T_{n}\left( \mathbf{p};\mathbf{a}%
,\left\vert \mathbf{b}\right\vert \right) \right\vert ,\left\vert
T_{n}\left( \mathbf{p};\left\vert \mathbf{a}\right\vert ,\left\vert \mathbf{b%
}\right\vert \right) \right\vert \right\} \geq 0.  \label{a.6}
\end{equation}

In 1990, Dragomir \cite{6b} considered the following class associated to a
pair of synchronous sequences $\mathbf{a},$ $\mathbf{b}$; 
\begin{equation*}
\bar{S}\left( \mathbf{a},\mathbf{b}\right) :=\left\{ \mathbf{x}\in \mathbb{R}%
^{n}|\left( \mathbf{a}+\mathbf{x},\mathbf{b}\right) \text{ and }\left( 
\mathbf{a}-\mathbf{x},\mathbf{b}\right) \text{ are synchronous}\right\} .
\end{equation*}%
It can be shown that $\bar{S}\left( \mathbf{a},\mathbf{b}\right) \neq
\emptyset $ and one has the representation 
\begin{equation}
T_{n}\left( \mathbf{p};\mathbf{a},\mathbf{b}\right) =\sup\limits_{\mathbf{x}%
\in \bar{S}\left( \mathbf{a},\mathbf{b}\right) }\left\vert T_{n}\left( 
\mathbf{p};\mathbf{x},\mathbf{b}\right) \right\vert \geq 0.  \label{a.7}
\end{equation}%
Now, if $\mathbf{k}=\left( k,k,\dots ,k\right) $ is a constant sequence and
if we denote by $\mathbf{a}\vee \mathbf{k:}=\left( \max \left\{
a_{1},k\right\} ,\dots ,\max \left\{ a_{n},k\right\} \right) $ and by $%
\mathbf{a}\wedge \mathbf{k:}=\left( \min \left\{ a_{1},k\right\} ,\dots
,\min \left\{ a_{n},k\right\} \right) ,$ then we may state the following
result obtained in the general setting of positive linear functionals by
Dragomir in 1993, \cite{8b} 
\begin{multline}
T_{n}\left( \mathbf{p};\mathbf{a},\mathbf{b}\right) \geq \max \left\{
\left\vert T_{n}\left( \mathbf{p};\mathbf{a}\vee \mathbf{k},\mathbf{b}%
\right) \right\vert +\left\vert T_{n}\left( \mathbf{p};\mathbf{a}\wedge 
\mathbf{k},\mathbf{b}\right) \right\vert \right. ,  \label{a.8} \\
\left. \left\vert T_{n}\left( \mathbf{p};\mathbf{a},\mathbf{b}\vee \mathbf{k}%
\right) \right\vert +\left\vert T_{n}\left( \mathbf{p};\mathbf{a},\mathbf{b}%
\wedge \mathbf{k}\right) \right\vert \right\} \geq 0,
\end{multline}%
provided $\mathbf{a}\ $and $\mathbf{b}$ are synchronous.

If $\mathbf{k}=\mathbf{0},$ and $\mathbf{a}_{+}:=\mathbf{a}\vee \mathbf{0},$ 
$\mathbf{a}_{-}:=\mathbf{a}\wedge \mathbf{0},$ then for synchronous
sequences $\mathbf{a},$ $\mathbf{b}$ one has 
\begin{multline}
T_{n}\left( \mathbf{p};\mathbf{a},\mathbf{b}\right) \geq \max \left\{
\left\vert T_{n}\left( \mathbf{p};\mathbf{a}_{+},\mathbf{b}\right)
\right\vert +\left\vert T_{n}\left( \mathbf{p};\mathbf{a}_{-},\mathbf{b}%
\right) \right\vert \right. ,  \label{a.9} \\
\left. \left\vert T_{n}\left( \mathbf{p};\mathbf{a},\mathbf{b}_{+}\right)
\right\vert +\left\vert T_{n}\left( \mathbf{p};\mathbf{a},\mathbf{b}%
_{-}\right) \right\vert \right\} \geq 0.
\end{multline}%
Note that, since, obviously 
\begin{align*}
\left\vert T_{n}\left( \mathbf{p};\mathbf{a}_{+},\mathbf{b}\right)
\right\vert +\left\vert T_{n}\left( \mathbf{p};\mathbf{a}_{-},\mathbf{b}%
\right) \right\vert & \geq \left\vert T_{n}\left( \mathbf{p};\mathbf{a}_{+},%
\mathbf{b}\right) +T_{n}\left( \mathbf{p};\mathbf{a}_{-},\mathbf{b}\right)
\right\vert  \\
& =\left\vert T_{n}\left( \mathbf{p};\left\vert \mathbf{a}\right\vert ,%
\mathbf{b}\right) \right\vert ,
\end{align*}%
then by (\ref{a.6}) and (\ref{a.9}), we deduce the sequence of inequalities 
\begin{equation}
T_{n}\left( \mathbf{p};\mathbf{a},\mathbf{b}\right) \geq \left\vert
T_{n}\left( \mathbf{p};\mathbf{a}_{+},\mathbf{b}\right) \right\vert
+\left\vert T_{n}\left( \mathbf{p};\mathbf{a}_{-},\mathbf{b}\right)
\right\vert \geq \left\vert T_{n}\left( \mathbf{p};\left\vert \mathbf{a}%
\right\vert ,\mathbf{b}\right) \right\vert \geq 0,  \label{a.10}
\end{equation}%
provided $\mathbf{a}$ and $\mathbf{b}$ are synchronous. This is a refinement
of (\ref{a.6}).

If one would like to drop the assumption of nonnegativity for the components
of $\mathbf{p},$ then one may state the following inequality obtained by
Mitrinovi\'{c} and Pe\v{c}ari\'{c} in 1991, \cite{MP}:

If $0\leq P_{i}\leq P_{n}$ for each $i\in \left\{ 1,\dots ,n-1\right\} ,$
then 
\begin{equation}
T_{n}\left( \mathbf{p};\mathbf{a},\mathbf{b}\right) \geq 0  \label{a.11}
\end{equation}%
provided $\mathbf{a}$ and $\mathbf{b}$ are sequences with the same
monotonicity.

If $\mathbf{a}$ and $\mathbf{b}$ are monotonic in the opposite sense, the
sign of the inequality (\ref{a.11}) reverses.

In this paper we point out other inequalities for the weighted \v{C}eby\v{s}%
ev's functional $T_{n}\left( \mathbf{p};\mathbf{a},\mathbf{b}\right) .$

\section{Some New Inequalities}

The following lemma holds.

\begin{lemma}
\label{l2.1}Let $\mathbf{a}=\left( a_{1},\dots ,a_{n}\right) ,$ $\mathbf{b}%
=\left( b_{1},\dots ,b_{n}\right) $ and $\mathbf{p}=\left( p_{1},\dots
,p_{n}\right) $ be sequences of real numbers. Define 
\begin{align*}
P_{i}& :=\sum_{k=1}^{i}p_{k},\;\;\bar{P}_{i}=P_{n}-P_{i}, \\
A_{i}\left( \mathbf{p}\right) & =\sum_{k=1}^{i}p_{k}a_{k},\;\;\bar{A}%
_{i}\left( \mathbf{p}\right) =A_{n}\left( \mathbf{p}\right) -A_{i}\left( 
\mathbf{p}\right) .
\end{align*}%
If we assume that $P_{n}\neq 0,$ then we have the identities 
\begin{align}
T_{n}\left( \mathbf{p};\mathbf{a},\mathbf{b}\right) & =\frac{1}{P_{n}^{2}}%
\sum_{i=1}^{n-1}\det \left[ 
\begin{array}{ll}
P_{i} & P_{n} \\ 
A_{i}\left( \mathbf{p}\right)  & A_{n}\left( \mathbf{p}\right) 
\end{array}%
\right] \cdot \Delta b_{i}  \label{2.1} \\
& =\frac{1}{P_{n}}\sum_{i=1}^{n-1}P_{i}\left[ \frac{A_{n}\left( \mathbf{p}%
\right) }{P_{n}}-\frac{A_{i}\left( \mathbf{p}\right) }{P_{i}}\right] \cdot
\Delta b_{i}  \notag \\
& =\frac{1}{P_{n}^{2}}\sum_{i=1}^{n-1}P_{i}\bar{P}_{i}\left[ \frac{\bar{A}%
_{i}\left( \mathbf{p}\right) }{\bar{P}_{i}}-\frac{A_{i}\left( \mathbf{p}%
\right) }{P_{i}}\right] \cdot \Delta b_{i},  \notag
\end{align}%
where $\Delta b_{i}:=b_{i+1}-b_{i}$ $\left( i=0,\dots ,n-1\right) $ is the
forward difference.
\end{lemma}

\begin{proof}
We use the following well known summation by parts formula 
\begin{equation}
\sum_{\ell =p}^{q-1}d_{\ell }\Delta v_{\ell }=d_{\ell }v_{\ell }\big|%
_{p}^{q}-\sum_{\ell =p}^{q-1}v_{\ell +1}\Delta d_{\ell },  \label{2.2}
\end{equation}
where $d_{\ell },v_{\ell }\in \mathbb{R}$, $\ell =p,\dots ,q$ ($q>p,$ $p,q$
are natural numbers).

If we choose in (\ref{2.2}), $p=1,q=n,$ $d_{i}=P_{i}A_{n}\left( \bar{p}%
\right) -P_{n}A_{i}\left( \bar{p}\right) $ and $v_{i}=b_{i}$ $\left(
i=1,\dots ,n\right) ,$ then we get 
\begin{align*}
& \sum_{i=1}^{n-1}\left[ P_{i}A_{n}\left( \mathbf{p}\right)
-P_{n}A_{i}\left( \mathbf{p}\right) \right] \cdot \Delta b_{i} \\
& =\left[ P_{i}A_{n}\left( \mathbf{p}\right) -P_{n}A_{i}\left( \mathbf{p}%
\right) \right] b_{i}\big|_{1}^{n}-\sum_{i=1}^{n-1}\Delta \left(
P_{i}A_{n}\left( \mathbf{p}\right) -P_{n}A_{i}\left( \mathbf{p}\right)
\right) b_{i+1} \\
& =\left[ P_{n}A_{n}\left( \mathbf{p}\right) -P_{n}A_{n}\left( \mathbf{p}%
\right) \right] b_{n}-\left[ P_{1}A_{n}\left( \mathbf{p}\right)
-P_{n}A_{1}\left( \mathbf{p}\right) \right] b_{1} \\
& \;\;\;\;\;\;\;\;\;\;\;\;-\sum_{i=1}^{n-1}\left[ P_{i+1}A_{n}\left( \mathbf{%
p}\right) -P_{n}A_{i+1}\left( \mathbf{p}\right) -P_{i}A_{n}\left( \mathbf{p}%
\right) +P_{n}A_{i}\left( \mathbf{p}\right) \right] b_{i+1} \\
& =P_{n}p_{1}a_{1}b_{1}-p_{1}b_{1}A_{n}\left( \mathbf{p}\right)
-\sum_{i=1}^{n-1}\left( p_{i+1}A_{n}\left( \mathbf{p}\right)
-P_{n}p_{i+1}a_{i+1}\right) b_{i+1}
\end{align*}%
\begin{align*}
& =P_{n}p_{1}a_{1}b_{1}-p_{1}b_{1}A_{n}\left( \mathbf{p}\right) -A_{n}\left( 
\mathbf{p}\right)
\sum_{i=1}^{n-1}p_{i+1}b_{i+1}+P_{n}\sum_{i=1}^{n-1}p_{i+1}a_{i+1}b_{i+1} \\
&
=P_{n}\sum_{i=1}^{n}p_{i}a_{i}b_{i}-\sum_{i=1}^{n}p_{i}a_{i}%
\sum_{i=1}^{n}p_{i}b_{i} \\
& =P_{n}^{2}T_{n}\left( \mathbf{p};\mathbf{a},\mathbf{b}\right) 
\end{align*}%
which produces the first identity in (\ref{2.1}).

The second and third are obvious and we omit the details.
\end{proof}

The following result holds.

\begin{theorem}
\label{t2.1}Let $\mathbf{a}=\left( a_{1},\dots ,a_{n}\right) ,$ $\mathbf{b}%
=\left( b_{1},\dots ,b_{n}\right) $ and $\mathbf{p}=\left( p_{1},\dots
,p_{n}\right) $ be sequences of real numbers. Assume that $p_{i}\geq 0$ $%
\left( i\in \left\{ 1,\dots ,n\right\} \right) $ such that $P_{i}\neq 0$ $%
\left( i\in \left\{ 1,\dots ,n\right\} \right) .$

If either

\begin{enumerate}
\item[$\left( i\right) $] $\mathbf{b}$ is increasing and $\mathbf{a}$ a
last-max in mean sequence, i.e., $\mathbf{a}$ \ satisfies the condition 
\begin{equation*}
\frac{A_{n}\left( \mathbf{p}\right) }{P_{n}}\geq \frac{A_{i}\left( \mathbf{p}%
\right) }{P_{i}}
\end{equation*}%
for each $i\in \left\{ 1,\dots ,n-1\right\} ;$

or

\item[$\left( ii\right) $] $\mathbf{b}$\hspace{0.05in}is decreasing and $%
\mathbf{a}$ is a first-max in mean sequence, i.e., 
\begin{equation*}
\frac{A_{n}\left( \mathbf{p}\right) }{P_{n}}\leq \frac{A_{i}\left( \mathbf{p}%
\right) }{P_{i}}
\end{equation*}%
for each $i\in \left\{ 1,\dots ,n-1\right\} ;$
\end{enumerate}

then one has the inequality 
\begin{eqnarray}
&&T_{n}\left( \mathbf{p};\mathbf{a},\mathbf{b}\right)   \label{2.1e} \\
&\geq &\max \left\{ \left\vert A_{n}\left( \mathbf{p};\mathbf{a},\mathbf{b}%
\right) \right\vert ,\left\vert A_{n}\left( \mathbf{p};\mathbf{a},\left\vert 
\mathbf{b}\right\vert \right) \right\vert ,\left\vert T_{n}\left( \mathbf{p};%
\mathbf{a},\left\vert \mathbf{b}\right\vert \right) \right\vert ,\right.  
\notag \\
&&\left. \left\vert D_{n}\left( \mathbf{p};\mathbf{a},\mathbf{b}\right)
\right\vert ,\left\vert D_{n}\left( \mathbf{p};\mathbf{a},\left\vert \mathbf{%
b}\right\vert \right) \right\vert \geq 0\right\} ;  \notag
\end{eqnarray}%
where 
\begin{equation*}
A_{n}\left( \mathbf{p};\mathbf{a},\mathbf{b}\right) =\frac{1}{P_{n}}%
\sum_{i=1}^{n-1}\left\vert A_{i}\left( \mathbf{p}\right) \right\vert \Delta
b_{i}-\frac{\left\vert A_{n}\left( \mathbf{p}\right) \right\vert }{P_{n}}%
\cdot \frac{1}{P_{n}}\sum_{i=1}^{n-1}P_{i}\Delta b_{i}
\end{equation*}%
and%
\begin{equation*}
D_{n}\left( \mathbf{p};\mathbf{a},\mathbf{b}\right) :=\frac{1}{P_{n}^{2}}%
\sum_{i=1}^{n-1}P_{i}\left\vert \bar{A}_{i}\left( \mathbf{p}\right)
\right\vert \Delta b_{i}-\frac{1}{P_{n}^{2}}\sum_{i=1}^{n-1}P_{i}\left\vert
A_{i}\left( \mathbf{p}\right) \right\vert \Delta b_{i}.
\end{equation*}
\end{theorem}

\begin{proof}
If either $\left( i\right) $ or $\left( ii\right) $ holds, then 
\begin{align*}
& \left( \frac{A_{n}\left( \mathbf{p}\right) }{P_{n}}-\frac{A_{i}\left( 
\mathbf{p}\right) }{P_{i}}\right) \left( b_{i+1}-b_{i}\right)  \\
& =\left\vert \left( \frac{A_{n}\left( \mathbf{p}\right) }{P_{n}}-\frac{%
A_{i}\left( \mathbf{p}\right) }{P_{i}}\right) \left( b_{i+1}-b_{i}\right)
\right\vert  \\
& \geq \left\{ 
\begin{array}{l}
\left\vert \left( \dfrac{\left\vert A_{n}\left( \mathbf{p}\right)
\right\vert }{P_{n}}-\dfrac{\left\vert A_{i}\left( \mathbf{p}\right)
\right\vert }{P_{i}}\right) \left( b_{i+1}-b_{i}\right) \right\vert  \\ 
\\ 
\left\vert \left( \dfrac{\left\vert A_{n}\left( \mathbf{p}\right)
\right\vert }{P_{n}}-\dfrac{\left\vert A_{i}\left( \mathbf{p}\right)
\right\vert }{P_{i}}\right) \left( \left\vert b_{i+1}\right\vert -\left\vert
b_{i}\right\vert \right) \right\vert  \\ 
\\ 
\left\vert \left( \dfrac{A_{n}\left( \mathbf{p}\right) }{P_{n}}-\dfrac{%
A_{i}\left( \mathbf{p}\right) }{P_{i}}\right) \left( \left\vert
b_{i+1}\right\vert -\left\vert b_{i}\right\vert \right) \right\vert 
\end{array}%
\right. 
\end{align*}%
for each $i\in \left\{ 1,\dots ,n-1\right\} .$

Multiplying by $P_{i}>0,$ summing over $i$ from $1$ to $n-1,$ and using the
generalised triangle inequality, we get 
\begin{align*}
T_{n}\left( \mathbf{p};\mathbf{a},\mathbf{b}\right) & =\frac{1}{P_{n}}%
\sum_{i=1}^{n-1}P_{i}\left\vert \left[ \left( \frac{A_{n}\left( \mathbf{p}%
\right) }{P_{n}}-\frac{A_{i}\left( \mathbf{p}\right) }{P_{i}}\right) \right]
\left( \Delta b_{i}\right) \right\vert  \\
& \geq \frac{1}{P_{n}}\times \left\{ 
\begin{array}{l}
\left\vert \dsum\limits_{i=1}^{n-1}P_{i}\left( \dfrac{\left\vert A_{n}\left( 
\mathbf{p}\right) \right\vert }{P_{n}}-\dfrac{\left\vert A_{i}\left( \mathbf{%
p}\right) \right\vert }{P_{i}}\right) \left( b_{i+1}-b_{i}\right)
\right\vert  \\ 
\\ 
\left\vert \dsum\limits_{i=1}^{n-1}P_{i}\left( \dfrac{\left\vert A_{n}\left( 
\mathbf{p}\right) \right\vert }{P_{n}}-\dfrac{\left\vert A_{i}\left( \mathbf{%
p}\right) \right\vert }{P_{i}}\right) \left( \left\vert b_{i+1}\right\vert
-\left\vert b_{i}\right\vert \right) \right\vert  \\ 
\\ 
\left\vert \dsum\limits_{i=1}^{n-1}P_{i}\left( \dfrac{A_{n}\left( \mathbf{p}%
\right) }{P_{n}}-\dfrac{A_{i}\left( \mathbf{p}\right) }{P_{i}}\right) \left(
\left\vert b_{i+1}\right\vert -\left\vert b_{i}\right\vert \right)
\right\vert 
\end{array}%
\right. 
\end{align*}%
from where we easily deduce the first three bounds in  (\ref{2.1e}). The
last two bounds may be obtained by utilising the second equality in (\ref%
{2.1}) and we omit the details.
\end{proof}

\begin{remark}
\label{r1}We observe that if $\mathbf{a}=\left( a_{1},\dots ,a_{n}\right) $
is monotonic increasing in mean for a given $\mathbf{p}$ positive, i.e., 
\begin{equation*}
\frac{1}{P_{i}}A_{i}\left( \mathbf{p}\right) \leq \frac{1}{P_{i+1}}%
A_{i+1}\left( \mathbf{p}\right) 
\end{equation*}%
then obviously 
\begin{equation}
\frac{1}{P_{i}}A_{i}\left( \mathbf{p}\right) \leq \frac{1}{P_{n}}A_{n}\left( 
\mathbf{p}\right)   \label{2.2a}
\end{equation}%
for each $i\in \left\{ 1,\dots ,n-1\right\} ,$ i.e., $\mathbf{a}$ is a
last-max in mean sequence for that specific weight vector $\mathbf{p.}$ The
converse is not true, generally.

We also note that if $\mathbf{a}$ is monotonic nondecreasing, then for any
positive $\mathbf{p}$, it is increasing in mean and, \textit{a fortiori, }a
last-max in mean sequence.
\end{remark}

\begin{remark}
\label{r2}We observe, for $\bar{A}_{i}\left( \mathbf{p}\right) :=A_{n}\left( 
\mathbf{p}\right) -A_{i}\left( \mathbf{p}\right) ,$ $i\in \left\{ 1,\dots
,n-1\right\} ,$ that 
\begin{equation*}
\frac{A_{n}\left( \mathbf{p}\right) }{P_{n}}-\frac{A_{i}\left( \mathbf{p}%
\right) }{P_{i}}=\frac{\bar{P}_{i}}{P_{n}}\left[ \frac{\bar{A}_{i}\left( 
\mathbf{p}\right) }{\bar{P}_{i}}-\frac{A_{i}\left( \mathbf{p}\right) }{P_{i}}%
\right] 
\end{equation*}%
for each $i\in \left\{ 1,\dots ,n-1\right\} ,$ and thus, if we assume that $%
\mathbf{p}$ is positive, then 
\begin{equation*}
\frac{A_{n}\left( \mathbf{p}\right) }{P_{n}}\geq \frac{A_{i}\left( \mathbf{p}%
\right) }{P_{i}}\text{ \hspace{0.05in}for every \hspace{0.05in}}i\in \left\{
1,\dots ,n-1\right\} 
\end{equation*}%
if and only if 
\begin{equation*}
\frac{\bar{A}_{i}\left( \mathbf{p}\right) }{\bar{P}_{i}}\geq \frac{%
A_{i}\left( \mathbf{p}\right) }{P_{i}}\text{ \hspace{0.05in}for every 
\hspace{0.05in}}i\in \left\{ 1,\dots ,n-1\right\} .
\end{equation*}
\end{remark}

If we would like to omit the assumption of positivity for the sequence $%
\mathbf{p}$, then the following result providing sufficient conditions for
the functional $T_{n}\left( \mathbf{p};\mathbf{a},\mathbf{b}\right) $ to be
positive (negative) holds.

\begin{theorem}
\label{t2.3}Let $\mathbf{a}=\left( a_{1},\dots ,a_{n}\right) ,$ $\mathbf{b}%
=\left( b_{1},\dots ,b_{n}\right) $ and $\mathbf{p}=\left( p_{1},\dots
,p_{n}\right) $ be $n-$tuples of real numbers. If $\mathbf{b}$ is monotonic
nondecreasing and either

\begin{enumerate}
\item[$\left( i\right) $] 
\begin{equation*}
\det \left( 
\begin{array}{ll}
P_{i} & P_{n} \\ 
A_{i}\left( \mathbf{p}\right)  & A_{n}\left( \mathbf{p}\right) 
\end{array}%
\right) \geq 0\text{ for each }i\in \left\{ 1,\dots ,n-1\right\} ;
\end{equation*}%
or

\item[$\left( ii\right) $] $P_{i}>0$ for any $i\in \left\{ 1,\dots
,n\right\} $ and $\mathbf{a}$ is a last-max in mean sequence 

or

\item[$\left( iii\right) $] $0<P_{i}<P_{n}$ for every $i\in \left\{ 1,\dots
,n-1\right\} $ and 
\begin{equation*}
\frac{\bar{A}_{i}\left( \mathbf{p}\right) }{P_{i}}\geq \frac{A_{i}\left( 
\mathbf{p}\right) }{P_{i}}\text{ \hspace{0.05in}for each \hspace{0.05in}}%
i\in \left\{ 1,\dots ,n-1\right\} ;
\end{equation*}%
then 
\begin{equation}
T_{n}\left( \mathbf{p};\mathbf{a},\mathbf{b}\right) \geq 0.  \label{2.4a}
\end{equation}
\end{enumerate}

If $\mathbf{b}$ is monotonic nonincreasing and either $\left( i\right) $ or $%
\left( ii\right) $ or $\left( iii\right) $ from above holds, then the
reverse inequality in (\ref{2.4a}) holds true.
\end{theorem}

The proof of the theorem follows from the identities incorporated in Lemma %
\ref{l2.1} and we omit the details.

\section{Some Inequalities for Convex (Concave) Sequences}

The following result holds.

\begin{theorem}
\label{t3.1}Let $\mathbf{a}=\left( a_{1},\dots ,a_{n}\right) $ and $\mathbf{b%
}=\left( b_{1},\dots ,b_{n}\right) $ be two sequences of real numbers and $%
\mathbf{p}=\left( p_{1},\dots ,p_{n}\right) $ a sequence of positive real
numbers.

If $\mathbf{b}$ is convex (concave), i.e., 
\begin{equation}
\frac{b_{i+2}+b_{i}}{2}\geq \left( \leq \right) b_{i+1}\;\;\;\text{for each 
\hspace{0.05in}}i\in \left\{ 1,\dots ,n-2\right\}   \label{3.1}
\end{equation}%
and $\mathbf{a}$ satisfies the property 
\begin{equation}
a_{i+1}\leq \left( \geq \right) \frac{A_{n}\left( \mathbf{p}\right) }{P_{n}}%
,\;\;\;\text{for each \hspace{0.05in}}i\in \left\{ 1,\dots ,n-1\right\} ;
\label{3.2}
\end{equation}%
then we have the inequality 
\begin{equation}
T_{n}\left( \mathbf{p};\mathbf{a},\mathbf{b}\right) \geq \frac{1}{\left(
n-1\right) }\left( b_{n}-b_{1}\right) \cdot \frac{1}{P_{n}}%
\dsum\limits_{i=1}^{n-1}\left( n-i\right) p_{i}\left[ \frac{A_{n}\left( 
\mathbf{p}\right) }{P_{n}}-a_{i}\right] .  \label{3.2.a}
\end{equation}
\end{theorem}

\begin{proof}
We know, by \v{C}eby\v{s}ev's inequality that if $\mathbf{\bar{z}}$ and $%
\mathbf{\bar{u}}$ have the same monotonicity, then 
\begin{equation}
\left( n-1\right) \dsum\limits_{i=1}^{n-1}z_{i}u_{i}\geq
\dsum\limits_{i=1}^{n-1}z_{i}\dsum\limits_{i=1}^{n-1}u_{i}.  \label{3.3}
\end{equation}%
Define $z_{i}:=b_{i+1}-b_{i}$ and $u_{i}:=P_{i}A_{n}\left( \mathbf{p}\right)
-A_{i}\left( \mathbf{p}\right) P_{n}$ for $i\in \left\{ 1,\dots ,n-1\right\}
.$ Then 
\begin{equation*}
z_{i+1}-z_{i}=2\left( \frac{b_{i+2}+b_{i}}{2}-b_{i+1}\right) \geq \left(
\leq \right) 0\text{ \hspace{0.05in}for each }i\in \left\{ 1,\dots
,n-1\right\} 
\end{equation*}%
and 
\begin{eqnarray*}
u_{i+1}-u_{i} &=&P_{i+1}A_{n}\left( \mathbf{p}\right) -A_{i+1}\left( \mathbf{%
p}\right) P_{n}-P_{i}A_{n}\left( \mathbf{p}\right) +A_{i}\left( \mathbf{p}%
\right) P_{n} \\
&=&p_{i+1}A_{n}\left( \mathbf{p}\right) -a_{i+1}p_{i+1}P_{n} \\
&=&p_{i+1}P_{n}\left( \frac{A_{n}\left( \mathbf{p}\right) }{P_{n}}%
-a_{i+1}\right)  \\
&\geq &\left( \leq \right) 0
\end{eqnarray*}%
for each \hspace{0.05in}$i\in \left\{ 1,\dots ,n-1\right\} ,$ showing that $%
\mathbf{\bar{z}}$ and $\mathbf{\bar{u}}$ have the same monotonicity.
Applying (\ref{3.3}) and the first identity in (\ref{2.1}), we have 
\begin{align*}
T_{n}\left( \mathbf{p};\mathbf{a},\mathbf{b}\right) & =\frac{1}{P_{n}^{2}}%
\dsum\limits_{i=1}^{n-1}\left( P_{i}A_{n}\left( \mathbf{p}\right)
-A_{i}\left( \mathbf{p}\right) P_{n}\right) \left( b_{i+1}-b_{i}\right)  \\
& \geq \frac{1}{\left( n-1\right) P_{n}^{2}}\dsum\limits_{i=1}^{n-1}\left(
P_{i}A_{n}\left( \mathbf{p}\right) -A_{i}\left( \mathbf{p}\right)
P_{n}\right) \dsum\limits_{i=1}^{n-1}\left( b_{i+1}-b_{i}\right) 
\end{align*}%
\begin{align*}
& =\frac{1}{\left( n-1\right) P_{n}^{2}}\left[ A_{n}\left( \mathbf{p}\right)
\dsum\limits_{i=1}^{n-1}P_{i}-P_{n}\dsum\limits_{i=1}^{n-1}A_{i}\left( 
\mathbf{p}\right) \right] \left( b_{n}-b_{1}\right)  \\
& =\frac{1}{\left( n-1\right) }\left( b_{n}-b_{1}\right) \left[ \frac{%
A_{n}\left( \mathbf{p}\right) }{P_{n}}\cdot \frac{1}{P_{n}}%
\dsum\limits_{i=1}^{n-1}\left( n-i\right) p_{i}-\frac{1}{P_{n}}%
\dsum\limits_{i=1}^{n-1}\left( n-i\right) p_{i}a_{i}\right]  \\
& =\frac{1}{\left( n-1\right) }\left( b_{n}-b_{1}\right) \frac{1}{P_{n}}%
\dsum\limits_{i=1}^{n-1}\left( n-i\right) p_{i}\left[ \frac{A_{n}\left( 
\mathbf{p}\right) }{P_{n}}-a_{i}\right] 
\end{align*}%
and the inequality (\ref{3.2.a}) is proved.
\end{proof}

The second result which does not require positivity for the weights $\mathbf{%
p}=\left( p_{1},\dots ,p_{n}\right) ,$ is enclosed in the following theorem.

\begin{theorem}
\label{t3.2}Let $\mathbf{a}$, $\mathbf{b}$ and $\mathbf{p}$ be sequences of
real numbers. Assume $P_{i}:=\sum_{k=1}^{i}p_{k}>0$ for $i=1,\dots ,n$, $%
\mathbf{b}$ is convex (concave) and $\mathbf{a}$ satisfies the following
monotonicity in mean condition 
\begin{equation}
\frac{A_{i}\left( \mathbf{p}\right) }{P_{i}}\geq \left( \leq \right) \frac{%
A_{i+1}\left( \mathbf{p}\right) }{\bar{P}_{i+1}},\;\;\;\text{for \hspace{%
0.05in}}i\in \left\{ 1,\dots ,n-1\right\} .  \label{3.4}
\end{equation}%
Then one has the inequality 
\begin{multline}
\quad T_{n}\left( \mathbf{p};\mathbf{a},\mathbf{b}\right)   \label{3.5} \\
\geq \frac{1}{\sum_{i=1}^{n-1}\left( n-i\right) p_{i}}\dsum%
\limits_{i=1}^{n-1}\left( n-i\right) p_{i}\left[ \frac{A_{n}\left( \mathbf{p}%
\right) }{P_{n}}-a_{i}\right] \cdot \left[ b_{n}-\frac{B_{n}\left( \mathbf{p}%
\right) }{P_{n}}\right] ,\quad 
\end{multline}%
where $B_{n}\left( \mathbf{p}\right) :=\sum_{i=1}^{n}p_{i}b_{i}.$
\end{theorem}

\begin{proof}
We use the following \v{C}eby\v{s}ev weighted inequality 
\begin{equation}
\dsum\limits_{i=1}^{n-1}q_{i}\dsum\limits_{i=1}^{n-1}q_{i}z_{i}u_{i}\geq
\dsum\limits_{i=1}^{n-1}q_{i}z_{i}\dsum\limits_{i=1}^{n-1}q_{i}u_{i},
\label{3.6}
\end{equation}
provided $q_{i}\geq 0$ and $\mathbf{\bar{z}}$, $\mathbf{\bar{u}}$ are
monotonic in the same sense.

Now, if we define $q_{i}:=P_{i},$ $z_{i}:=b_{i+1}-b_{i}$ and $u_{i}:=\frac{%
A_{n}\left( \mathbf{p}\right) }{P_{n}}-\frac{A_{i+1}\left( \mathbf{p}\right) 
}{\bar{P}_{i+1}}$ for $i\in \left\{ 1,\dots ,n-1\right\} ,$ then by \v{C}eby%
\v{s}ev's inequality, (\ref{3.6}) and the second identity in (\ref{2.1}), we
have 
\begin{align}
T_{n}\left( \mathbf{p};\mathbf{a},\mathbf{b}\right) & =\frac{1}{P_{n}}%
\dsum\limits_{i=1}^{n-1}P_{i}\left[ \frac{A_{n}\left( \mathbf{p}\right) }{%
P_{n}}-\frac{A_{i}\left( \mathbf{p}\right) }{P_{i}}\right] \Delta b_{i}
\label{3.7} \\
& \geq \frac{1}{P_{n}\sum_{i=1}^{n-1}P_{i}}\dsum\limits_{i=1}^{n-1}P_{i}%
\left[ \frac{A_{n}\left( \mathbf{p}\right) }{P_{n}}-\frac{A_{i}\left( 
\mathbf{p}\right) }{P_{i}}\right] \dsum\limits_{i=1}^{n-1}P_{i}\Delta b_{i}.
\notag
\end{align}%
Since 
\begin{equation*}
\dsum\limits_{i=1}^{n-1}P_{i}=\dsum\limits_{i=1}^{n-1}\left( n-i\right)
p_{i},
\end{equation*}%
\begin{equation*}
\dsum\limits_{i=1}^{n-1}P_{i}\left[ \frac{A_{n}\left( \mathbf{p}\right) }{%
P_{n}}-\frac{A_{i}\left( \mathbf{p}\right) }{P_{i}}\right]
=\dsum\limits_{i=1}^{n-1}\left( n-i\right) p_{i}\left[ \frac{A_{n}\left( 
\mathbf{p}\right) }{P_{n}}-a_{i}\right] 
\end{equation*}%
and 
\begin{align*}
\dsum\limits_{i=1}^{n-1}P_{i}\Delta b_{i}& =P_{i}b_{i}\bigg|%
_{1}^{n}-\dsum\limits_{i=1}^{n-1}b_{i+1}\Delta P_{i} \\
& =P_{n}b_{n}-\dsum\limits_{i=1}^{n}b_{i}p_{i}
\end{align*}%
thus, by (\ref{3.7}), we get 
\begin{align*}
& T_{n}\left( \mathbf{p};\mathbf{a},\mathbf{b}\right)  \\
& \geq \frac{1}{P_{n}\sum_{i=1}^{n-1}\left( n-i\right) p_{i}}%
\dsum\limits_{i=1}^{n-1}\left( n-i\right) p_{i}\left[ \frac{A_{n}\left( 
\mathbf{p}\right) }{P_{n}}-a_{i}\right] \cdot \left[ P_{n}b_{n}-\dsum%
\limits_{i=1}^{n}b_{i}p_{i}\right]  \\
& =\frac{1}{\sum_{i=1}^{n-1}\left( n-i\right) p_{i}}\dsum\limits_{i=1}^{n-1}%
\left( n-i\right) p_{i}\left[ \frac{A_{n}\left( \mathbf{p}\right) }{P_{n}}%
-a_{i}\right] \cdot \left[ b_{n}-\frac{1}{P_{n}}\dsum%
\limits_{i=1}^{n}p_{i}b_{i}\right] 
\end{align*}%
and the theorem is completely proved.
\end{proof}

\begin{acknowledgement}
The author would like to thank the anonymous referee for valuable comments
that have been incorporated in the final version of the paper.
\end{acknowledgement}

\end{document}